\documentclass{amsart}
\usepackage{amsmath,amsfonts,amssymb}
\usepackage{amsthm}
\usepackage{hyperref}

\def\R{{\mathbb R}}

\def\virgp{\raise 2pt\hbox{,}}

\def\({\left(}
\def\){\right)}
\def\<{\left\langle}
\def\>{\right\rangle}
\def\le{\leqslant}
\def\ge{\geqslant}

\def\Eq#1#2{\mathop{\sim}\limits_{#1\rightarrow#2}}
\def\Tend#1#2{\mathop{\longrightarrow}\limits_{#1\rightarrow#2}}

\def\si{{\sigma}}

\def\F{\mathcal F}

\DeclareMathOperator{\RE}{Re}

\theoremstyle{plain}
\newtheorem{theorem}{Theorem}[section]
\newtheorem{lemma}[theorem]{Lemma}

\theoremstyle{definition}

\theoremstyle{remark}
\newtheorem{remark}[theorem]{Remark}
\newtheorem*{remark*}{Remark}

\numberwithin{equation}{section}

\begin{document}

\title[A nonlinear Poisson formula for the Schr\"odinger operator]{A
  nonlinear Poisson formula for the Schr\"odinger operator}    
\author[R. Carles]{R{\'e}mi Carles}
\address{CNRS \& Universit\'e Montpellier~2\\Math\'ematiques,
  CC~051\\Place Eug\`ene 
  Bataillon\\34095 
  Montpellier cedex 5, France}
\email{Remi.Carles@math.cnrs.fr}
\author[T. Ozawa]{Tohru Ozawa}
\address{Department of Mathematics\\ Hokkaido University\\ Sapporo
  060-810, Japan}
\email{ozawa@math.sci.hokudai.ac.jp}
\begin{abstract}
We prove a nonlinear Poisson type formula for the Schr\"odinger
group. Such a formula had been derived in a previous paper by the
authors, as a consequence of the study of the asymptotic behavior of
nonlinear wave operators 
for small data. In this note, we propose a direct proof, and extend
the range allowed for the power of the nonlinearity to the set of
all short range nonlinearities. Moreover, $H^1$-critical
nonlinearities are allowed.  
\end{abstract}
\maketitle

\section{Introduction}
For $n\ge 1$, define the Schr\"odinger group as $U(t)=
e^{i\frac{t}{2}\Delta}$, where $\Delta$ stands for the Laplacian of
$\R^n$. We normalize the Fourier transform on $\R^n$ as follows:
\begin{equation}\label{eq:fourier}
  \F f(\xi)=\widehat
  f(\xi)=\frac{1}{(2\pi)^{n/2}}\int_{\R^n}f(x)e^{-ix\cdot \xi} dx.
\end{equation}
For $r\ge 2$, we define
\begin{equation*}
  \delta(r) =
\frac{n}{2}-\frac{n}{r}\cdot
\end{equation*}
The main result of this note is:
\begin{theorem}\label{theo:main}
  Let $n\ge 1$, and fix $1+2/n<p<\infty$ if $n\le 2$, $1+2/n<p\le
  1+4/(n-2)$ if $n\ge 3$.  
Then for every $\phi \in X_p$, and almost
  all $\xi\in \R^n$, the following identity holds:
  \begin{equation}\label{eq:poisson}
    \begin{aligned}
      \int_0^{\pm \infty}e^{i\frac{t}{2}|\xi|^2}& \F
    \(|U(t)\phi|^{p-1}U(t)\phi \)(\xi)dt =\\
=& \int_0^{\pm
    \infty}|t|^{n(p-1)/2-2} U(t)\( \left|U(-t)\widehat
    \phi\right|^{p-1}U(-t)\widehat \phi \)(\xi)dt,
    \end{aligned}
    \end{equation}
where the space $X_p$ is defined as follows:
\begin{itemize}
\item If $1+2/n<p< 1+4/n$, then $X_p = \{ f\in L^2(\R^n)\ ;\
  |x|^{\delta(2p)}f\in L^2(\R^n)\}$.
\item If $p=1+4/n$, then $X_p =L^2(\R^n)$. 
\item If $p> 1+4/n$, then $X_p = H^{\delta(p+1)}(\R^n)$, the
  inhomogeneous   Sobolev space. 
\end{itemize}
\end{theorem}
For $p=1+4/n$, the above result was proved in \cite{COMRL}. It
was also established for $1+2/n<p<1+4/n$ if $n\le 2$, and
$1+4/(n+2)<p<1+4/n$ if 
$n\ge 3$, provided that $\phi\in H^1\cap \F(H^1)$. The proof in
\cite{COMRL} relies on pseudo-conformal invariances for the nonlinear
Schr\"odinger equation, as well as the explicit computation of the
first non-trivial term in the asymptotic expansion of nonlinear wave
operators near the origin. In this note, we provide a direct proof of
the above identity, which relies on the usual factorization of the
Schr\"odinger group. Moreover, we extend the range of values allowed
for $p$, and we consider a broader class (when $p\not =1+4/n$) for the
function $\phi$. We also show that both terms in \eqref{eq:poisson}
become infinite when $p=1+2/n$ and $\phi$ is a Gaussian function (see
\S\ref{sec:long}). Note 
that $p=1+2/n$ corresponds to the long range case for the scattering
theory associated to the nonlinear Schr\"odinger equation with
nonlinearity $|u|^{p-1}u$; see e.g. \cite{CazCourant,Ginibre} and
references therein. 
\bigbreak

Let us point out some similarities between \eqref{eq:poisson} and the
usual Poisson formula. First, if we write $U(t) = \F^{-1}
e^{-i\frac{t}{2}|\xi|^2}\F$, we see that the right hand side of
\eqref{eq:poisson} has an additional Fourier transform compared to the
left hand side. Moreover, we will see in the proof that
\eqref{eq:poisson} relies on an inversion $t\mapsto 1/t$. This is the
same as for the Poisson formula associated to the heat equation, or to
the Jacobi theta function; see
e.g. \cite{Tata1}. 
\smallbreak

Note also that the definition of the space $X_p$ (which will become
natural in the course of the proof of the above result) is reminiscent
of the discussion related to scattering theory for the nonlinear
Schr\"odinger equation with nonlinearity $|u|^{p-1}u$. The case
$p=1+4/n$ corresponds to the $L^2$-critical nonlinearity. For
$p<1+4/n$, it is usual to work in weighted $L^2$ spaces, while for
$p>1+4/n$, Sobolev spaces are more convenient (see
e.g. \cite{NakanishiOzawa}). Also, note that for $p>1+4/n$, the
upper bound for $p$ allows $H^1$-critical nonlinearities
($p=1+4/(n-2)$ for $n\ge 3$), thanks to endpoint Strichartz
estimates. As mentioned above, when $p$ reaches the 
long range case $p=1+2/n$, \eqref{eq:poisson} becomes irrelevant.

\section{Proof of Theorem~\ref{theo:main}}

We recall the classical factorization of the Schr\"odinger group:
$U(t)= M_tD_t\F M_t$, where $M_t$ is the multiplication by
$e^{i|x|^2/(2t)}$, $\F$ is the Fourier transform \eqref{eq:fourier},
and $D_t$ is the dilation operator
\begin{equation*}
  D_t f(x)= \frac{1}{(it)^{n/2}}f\(\frac{x}{t}\). 
\end{equation*}
We first prove that both terms in \eqref{eq:poisson} are well defined
for $\phi\in X_p$ and almost all $\xi\in \R^n$:
\begin{lemma}
  Let $p$ as in Theorem~\ref{theo:main}, and $\phi\in X_p$. Let $F$
  denote either of the two 
  terms in \eqref{eq:poisson}. Then $F\in L^2(\R^n)$. More precisely,
  there exists $C>0$ independent of $\phi\in X_p$ such that:
  \begin{equation*}
    \|F\|_{L^2}\le C
\left\{
  \begin{aligned}
    &\left\lVert |x|^{\delta(2p)}\phi\right\rVert_{L^2}^{\theta
    p}\|\phi\|_{L^2}^{(1-\theta) p} && \text{if }1+2/n<p<1+4/n,\\
&\|\phi\|_{L^2}&& \text{if } p=1+4/n,\\
&\left\lVert (-\Delta)^{\delta(p+1)/2}\phi\right\rVert_{L^2}^{(1-\si)
    p}\|\phi\|_{L^2}^{\si p} && \text{if }p>1+4/n,
  \end{aligned}
\right.
  \end{equation*}
where $\theta$ and $\si$ are given by:
\begin{equation*}
  \theta = \frac{4}{n(p-1)}-1\quad ;\quad \si =
  \frac{n+4-(n-4)p}{np(p-1)}.  
\end{equation*}
\end{lemma}
\begin{remark}
  We check the following algebraic identities:
  \begin{itemize}
  \item $0<\theta<1 \iff 1+2/n<p<1+4/n$, and $\theta =0 \iff p=1+4/n$.
\item $\si<1 \iff p>1+4+n$; $\si=1\iff p=1+4+n$.
\item $\si>0$, since for $n\ge 3$,
  $(n-2)p\le n+2$.  
  \end{itemize}
\end{remark}
\begin{proof} 
By symmetry, we consider only the plus sign in
  \eqref{eq:poisson}. We distinguish three cases, according to the
  value of $p$. 

\noindent {\bf First case: $1+2/n<p<1+4/n$.} Let $\psi \in L^2(\R^n)$,
$T>0$, and $q$ be defined by $2/q=\delta(p+1)$. Note that $0<2/q<1$,
so the pair $(q,p+1)$ is Strichartz admissible. By duality, we have:
\begin{align*}
  \Big\lvert \Big\langle\int_0^T e^{i\frac{t}{2}|\xi|^2}\F
    \big(|U(t)\phi|^{p-1}&U(t)\phi \big)dt, \widehat
    \psi\Big\rangle\Big\rvert=  
\Big\lvert \int_0^T \<
    |U(t)\phi|^{p-1}U(t)\phi , U(t) \psi\>dt\Big\rvert\\
&\le \int_0^T \left\lVert U(t)\phi\right\rVert_{L^{p+1}}^p \left\lVert
    U(t)\psi\right\rVert_{L^{p+1}}dt\\
&\le T^{1-n(p-1)/4}\left\lVert
    U(\cdot)\phi\right\rVert_{L^q(0,T;L^{p+1})}^p \left\lVert
    U(\cdot)\psi\right\rVert_{L^q(0,T;L^{p+1})}\\
&\le C  T^{1-n(p-1)/4}\left\lVert
    \phi\right\rVert_{L^2}^p \left\lVert
    \psi\right\rVert_{L^2},
\end{align*}
where $C$, independent of $T$, is provided by Strichartz
inequalities. Note that for the H\"older inequality in time, we have
used the formula:
\begin{equation*}
  1= \(1 -\frac{n(p-1)}{4}\) + \frac{p+1}{q}.
\end{equation*}
We also have directly
\begin{align*}
  \left\lVert \int_T^\infty e^{i\frac{t}{2}|\xi|^2}\F
    \big(|U(t)\phi|^{p-1}U(t)\phi \big)dt \right\rVert_{L^2}&\le 
\int_T^\infty\left\lVert e^{i\frac{t}{2}|\xi|^2}\F 
    \big(|U(t)\phi|^{p-1}U(t)\phi \big)\right\rVert_{L^2} dt \\
&\le \int_T^\infty \left\lVert U(t)\phi\right\rVert_{L^{2p}}^p dt\\
&\le C \int_T^\infty  \left\lVert
    U(t)\phi\right\rVert_{L^{2p}}^p dt.
\end{align*}
Using the factorization for the group $U$ recalled above, we find:
\begin{align*}
  \left\lVert
    U(t)\phi\right\rVert_{L^{2p}} = t^{-\delta(2p)}\left\lVert
    \F M_t\phi\right\rVert_{L^{2p}}\le C t^{-\delta(2p)} \left\lVert
    \F M_t\phi\right\rVert_{\dot H^{\delta(2p)}},  
\end{align*}
where we have used the critical Sobolev embedding. We infer
\begin{align*}
  \left\lVert \int_T^\infty e^{i\frac{t}{2}|\xi|^2}\F
    \big(|U(t)\phi|^{p-1}U(t)\phi \big)dt \right\rVert_{L^2}&\le 
C\int_T^\infty  t^{-p\delta(2p)} \left\lVert
    \lvert x\rvert^{\delta(2p)}\phi\right\rVert_{L^2}^p dt\\
&\le CT^{1-n(p-1)/2} \left\lVert
    \lvert x\rvert^{\delta(2p)}\phi\right\rVert_{L^2}^p.
\end{align*}
We have finally, for any $T>0$:
\begin{align*}
 \left\lVert \int_0^\infty e^{i\frac{t}{2}|\xi|^2}\F
    \big(|U(t)\phi|^{p-1}U(t)\phi \big)dt \right\rVert_{L^2} \le
    C\Big(&T^{1-n(p-1)/4}\left\lVert 
    \phi\right\rVert_{L^2}^p\\
& + T^{1-n(p-1)/2} \left\lVert
    \lvert x\rvert^{\delta(2p)}\phi\right\rVert_{L^2}^p\Big),
\end{align*}
where $C$ is independent of $T$. Optimizing in $T$, we find:
\begin{align*}
 \left\lVert \int_0^\infty e^{i\frac{t}{2}|\xi|^2}\F
    \big(|U(t)\phi|^{p-1}U(t)\phi \big)dt \right\rVert_{L^2} \le
    C\left\lVert
    \lvert x\rvert^{\delta(2p)}\phi\right\rVert_{L^2}^{\theta p}\left\lVert 
    \phi\right\rVert_{L^2}^{(1-\theta)p} ,
\end{align*}
where $\theta= \frac{4}{n(p-1)}-1$. For the other term involved in
\eqref{eq:poisson}, we proceed in a similar fashion:
\begin{align*}
   &\Big\lvert \Big\langle\int_{1/T}^\infty t^{n(p-1)/2-2}
    U(t)\Big( \left|U(-t)\widehat
    \phi\right|^{p-1}  U(-t)\widehat \phi \Big)dt, \widehat
    \psi\Big\rangle\Big\rvert=  \\
=&\Big\lvert \int_{1/T}^\infty t^{n(p-1)/2-2} \Big\langle 
    \left|U(-t)\widehat
    \phi\right|^{p-1}  U(-t)\widehat \phi , U(-t)\widehat
    \psi\Big\rangle dt\Big\rvert\\
\le & \int_{1/T}^\infty t^{n(p-1)/2-2} \left\lVert U(-t)\widehat \phi
    \right\rVert_{L^{p+1}}^p \left\lVert U(-t)\widehat \psi
    \right\rVert_{L^{p+1}}dt \\
\le & \(\int_{1/T}^\infty
    t^{\(n(p-1)/2-2\)/\(1-n(p-1)/4\)}dt\)^{1-n(p-1)/4}
\left\lVert
    U(\cdot)^{-1}\widehat \phi\right\rVert_{L^qL^{p+1}}^p
    \left\lVert 
    U(\cdot)^{-1}\widehat \psi\right\rVert_{L^q L^{p+1}}\\
\le & \ C T^{1-n(p-1)/4}\left\lVert
     \phi\right\rVert_{L^2}^p
    \left\lVert \psi\right\rVert_{L^2}, 
\end{align*}
for the same $q$ as above, given by $2/q=\delta(p+1)$. We also have
directly
\begin{align*}
  &\Big\lVert \int_0^{1/T}t^{n(p-1)/2-2}
    U(t)\Big( \left|U(-t)\widehat
    \phi\right|^{p-1}  U(-t)\widehat \phi \Big)dt\Big\rVert_{L^2}\le \\
\le & \int_0^{1/T}t^{n(p-1)/2-2} \left\lVert\left|U(-t)\widehat
    \phi\right|^{p-1}  U(-t)\widehat \phi\right\rVert_{L^2} dt\\
\le&   
\int_0^{1/T}t^{n(p-1)/2-2} \left\lVert  U(-t)\widehat
    \phi\right\rVert_{L^{2p}}^p dt \\
\le&   
C \int_0^{1/T}t^{n(p-1)/2-2} \left\lVert  U(-t)\widehat
    \phi\right\rVert_{\dot H^{\delta(2p)}}^p dt\\
\le&   
C\int_0^{1/T}t^{n(p-1)/2-2} \left\lVert  \widehat
    \phi\right\rVert_{\dot H^{\delta(2p)}}^p dt =
    CT^{1-n(p-1)/2}\left\lVert  \lvert x\rvert^{\delta(2p)}
    \phi\right\rVert_{L^2}^p.
\end{align*}
We infer:
\begin{align*}
  \Big\lVert \int_0^\infty t^{n(p-1)/2-2}
    U(t)\Big( \left|U(-t)\widehat
    \phi\right|^{p-1}  U(-t)\widehat \phi \Big)dt\Big\rVert_{L^2}&\le
    C\Big( T^{1-n(p-1)/4}\left\lVert
     \phi\right\rVert_{L^2}^p\\
+& T^{1-n(p-1)/2}\left\lVert  \lvert x\rvert^{\delta(2p)}
    \phi\right\rVert_{L^2}^p\Big). 
\end{align*}
We can then conclude as above.

\noindent {\bf Second case: $p=1+4/n$.} In this case, note that the
power of $t$ in the second term of \eqref{eq:poisson} is zero:
$n(p-1)/2-2=0$. To prove the result in this case, just notice that the
above proof remains valid: for $\psi\in L^2(\R^n)$ and $T>0$, we now
have
\begin{align*}
  \Big\lvert \Big\langle\int_0^T e^{i\frac{t}{2}|\xi|^2}\F
    \big(|U(t)\phi|^{p-1}U(t)\phi \big)dt, \widehat
    \psi\Big\rangle\Big\rvert &\le C  T^{1-n(p-1)/4}\left\lVert
    \phi\right\rVert_{L^2}^p \left\lVert
    \psi\right\rVert_{L^2}\\
&\le C  \left\lVert
    \phi\right\rVert_{L^2}^p \left\lVert
    \psi\right\rVert_{L^2},
\end{align*}
where $C$ is independent of $T$. The estimate for the
other term in \eqref{eq:poisson} is straightforward, by duality.

\noindent {\bf Third case: $p>1+4/n$.} For $\psi\in L^2(\R^n)$, we
compute 
\begin{align*}
  \Big\lvert \Big\langle\int_0^\infty e^{i\frac{t}{2}|\xi|^2}\F
    \big(|U(t)\phi|^{p-1} U(t)\phi \big)dt, \widehat
    \psi\Big\rangle\Big\rvert &\le \int_0^\infty \left\lVert
    U(t)\phi\right\rVert_{L^{p+1}}^p \left\lVert 
    U(t)\psi\right\rVert_{L^{p+1}}dt\\
\le \left\lVert
    U(\cdot)\phi\right\rVert_{L^\infty L^{p+1}}^{(1-\si) p}&\left\lVert
    U(\cdot)\phi\right\rVert_{L^q L^{p+1}}^{\si p} \left\lVert 
    U(\cdot)\psi\right\rVert_{L^qL^{p+1}},
\end{align*}
for $2/q= \delta(p+1)$, 
where we have used the identity $1=\si p/q + 1/q$. We conclude thanks
to the Sobolev embedding $\dot H^{\delta(p+1)}\hookrightarrow L^{p+1}$
and Strichartz inequalities. Note that for $n\ge 3$ and $p=1+4/(n-2)$,
we use endpoint estimates \cite{KT}. 

For the other term, write
\begin{align*}
  \int_0^\infty & t^{n(p-1)/2-2}\left\lVert
    U(-t)\widehat \phi\right\rVert_{L^{p+1}}^p \left\lVert 
    U(-t)\widehat \psi\right\rVert_{L^{p+1}}dt\le \\
&\le \(\sup_{t>0}t^{n(p-1)/2-2}\left\lVert
    U(-t)\widehat \phi\right\rVert_{L^{p+1}}^{(1-\si)p}\) \left\lVert
    U(\cdot)^{-1}\widehat \phi\right\rVert_{L^q L^{p+1}}^{\si p} \left\lVert 
    U(\cdot)^{-1}\widehat \psi\right\rVert_{L^q L^{p+1}} \\
& \le C\(\sup_{t>0}t^{n(p-1)/2-2}\left\lVert
    U(-t)\widehat \phi\right\rVert_{L^{p+1}}^{(1-\si)p}\) \lVert
    \phi\rVert_{L^2}^{\si p} \lVert  \psi\rVert_{L^2}.
\end{align*}
We then remark that
\begin{align*}
 \left\lVert
    U(-t)\widehat \phi\right\rVert_{L^{p+1}}&= \left\lVert
    M_t U(-t)\widehat \phi\right\rVert_{L^{p+1}}=
\left\lVert
    D_{-t}\F M_{-t}\widehat \phi\right\rVert_{L^{p+1}}=
    \frac{1}{|t|^{\delta(p+1)}}\left\lVert 
    \F M_{-t}\widehat \phi\right\rVert_{L^{p+1}}  \\
& \le \frac{C}{|t|^{\delta(p+1)}} \left\lVert
    \F M_{-t}\widehat \phi\right\rVert_{\dot H^{\delta(p+1)}}
= \frac{C}{|t|^{\delta(p+1)}} \left\lVert
    \lvert x\rvert^{\delta(p+1)} M_{-t}\widehat
    \phi\right\rVert_{L^2}\\
&\le  \frac{C}{|t|^{\delta(p+1)}} \left\lVert
    \lvert x\rvert^{\delta(p+1)}\widehat
    \phi\right\rVert_{L^2}=\frac{C}{|t|^{\delta(p+1)}} \left\lVert
    \phi\right\rVert_{\dot H^{\delta(p+1)}} .
\end{align*}
In view of the identity $(1-\si)p\delta(p+1)= n(p-1)/2-2$, this yields
\begin{equation*}
 \left\lVert \int_0^{\pm
    \infty}|t|^{n(p-1)/2-2} U(t)\( \left|U(-t)\widehat
    \phi\right|^{p-1}U(-t)\widehat \phi \)dt\right\rVert_{L^2} \le C   
\left\lVert \phi\right\rVert_{\dot H^{\delta(p+1)}}^{(1-\si)
    p}\|\phi\|_{L^2}^{\si p},
\end{equation*}
which completes the proof of the lemma. 
\end{proof}
We can now prove Theorem~\ref{theo:main}.
\begin{proof}[Proof of Theorem~\ref{theo:main}]
  Recall the decomposition $U(t)=M_t D_t \F M_t$. Direct computations
  yield:
  \begin{align}
    \F D_t &=D_{1/t}\F,\label{eq:3}\\
D_t^{-1}&= i^n D_{1/t},\label{eq:4}\\
\F^{-1}D_t^{-1}&= i^n D_t \F^{-1}. \label{eq:5}
  \end{align}
We infer
\begin{equation*}
  U(-t)= U(t)^{-1}= M_{-t}\F^{-1}D_t^{-1}M_{-t} = i^n M_{-t}D_t
  \F^{-1}M_{-t}. 
\end{equation*}
Since $U(t)= \F^{-1}M_{-1/t}\F$, we deduce
\begin{equation*}
 U(-t)\F = i^n  M_{-t}D_t U\(\frac{1}{t}\),
\end{equation*}
which in turn implies
\begin{align*}
  \left\lvert U(-t)\widehat
    \phi\right\rvert^{p-1}U(-t)\widehat \phi &= i^n M_{-t} \left\lvert
    D_t U\(\frac{1}{t}\)\phi\right\rvert^{p-1}D_t
    U\(\frac{1}{t}\)\phi\\
&= i^n t^{-n(p-1)/2}M_{-t}D_t \(\left\lvert
    U\(\frac{1}{t}\)\phi\right\rvert^{p-1}
    U\(\frac{1}{t}\)\phi\) .  
\end{align*}
Using \eqref{eq:3} and \eqref{eq:4} again, we have then:
\begin{equation}\label{eq:poissonponct}
  U(t)\(\left\lvert U(-t)\widehat
    \phi\right\rvert^{p-1}U(-t)\widehat \phi\) = t^{-n(p-1)/2}M_t
    \F\(\left\lvert U\(\frac{1}{t}\)\phi\right\rvert^{p-1}
    U\(\frac{1}{t}\)\phi\). 
\end{equation}
Theorem~\ref{theo:main} follows by integrating the above identity on a
half line, and using the change of variable $t\mapsto 1/t$.
\end{proof}
\begin{remark}
  The identity \eqref{eq:poissonponct} can also be considered as a
  Poisson formula, by writing $U(t)$ on the left hand side, and $\F$
  on the right hand side, as integrals. 
\end{remark}

\section{The long range case}
\label{sec:long}

When $\phi$ is a Gaussian function, the value in \eqref{eq:poisson}
can be computed explicitly. 
For $\RE z>0$,  define:
\begin{equation*}
  g_z (x) = e^{-z\frac{|x|^2}{2}},\quad x\in \R^n.
\end{equation*}
We have:
\begin{equation*}
  \int_{\R^n}g_z (x) dx = \( \frac{2\pi}{z}\)^{n/2}.
\end{equation*}
We compute:
\begin{equation*}
  \F g_z (\xi) = z^{-n/2} e^{-\frac{|\xi|^2}{2z}}, 
\end{equation*}
and
\begin{equation*}
  U(t)g_z(x) = \(1+itz\)^{-n/2} e^{-\frac{z}{1+itz}\frac{|x|^2}{2}}. 
\end{equation*}
Note that if $z=a+ib$,
\begin{equation*}
  \RE \( \frac{z}{1+itz}\) = \frac{a}{(1-tb)^2 +a^2t^2}>0.
\end{equation*}
For $p>1$ and $z=a+ib$, we find:
\begin{equation*}
  \left\lvert U(t)g_z \right\rvert^{p-1}U(t)g_z  = \frac{e^{-\frac{(p-1)
  a}{(1-bt)^2+(at)^2}\frac{|x|^2}{2}}}{\(
  (1-bt)^2+(at)^2\)^{n(p-1)/4}}  (1+itz)^{-n/2}
  e^{-\frac{z}{1+itz}\frac{|x|^2}{2}} .
\end{equation*}
Set 
\begin{equation*}
  \zeta = \frac{(p-1)
  a}{(1-bt)^2+(at)^2}+ \frac{z}{1+itz}. 
\end{equation*}
We have:
\begin{equation*}
  \F\( \left\lvert U(t)g_z \right\rvert^{p-1}U(t)g_z\) = \frac{1}{\(
  (1-bt)^2+(at)^2\)^{n(p-1)/4}} (1+itz)^{-n/2} \zeta^{-n/2}
  e^{-\frac{|x|^2}{2\zeta}}.  
\end{equation*}
Consider the case $z\in \R$: $b=0$. We find:
\begin{equation*}
  \zeta = \frac{a}{1+(at)^2}\(p -iat\).
\end{equation*}
We infer:
\begin{equation*}
  e^{i\frac{t}{2}|x|^2}\F\( \left\lvert U(t)g_z
\right\rvert^{p-1}U(t)g_z\) = \frac{1}{
  \(1+(at)^2\)^{n(p-1)/4}}\(\zeta
(1+ita)\)^{-n/2}e^{\(it-\frac{1}{\zeta}\)\frac{|x|^2}{2}}. 
\end{equation*}
We compute
\begin{equation*}
  it-\frac{1}{\zeta} =
  \frac{iat p-1}{a(p-iat)}\Tend t \infty
  - \frac{p}{a}. 
\end{equation*}
Also,
\begin{equation*}
  \zeta 
(1+ita) = \frac{a}{1+(at)^2} (p -iat) (1+ita)\Tend t\infty a. 
\end{equation*}
We have finally:
\begin{equation*}
  e^{i\frac{t}{2}|x|^2}\F\( \left\lvert U(t)g_z
\right\rvert^{2\si}U(t)g_z\)\Eq t\infty
\frac{1}{\(1+(at)^2\)^{n(p-1)/4}}\frac{1}{a^{n/2}}
e^{-p\frac{|x|^2}{2a}}. 
\end{equation*}
Integrating with respect to $t$, the integral is convergent if and only
if $p>1+2/n$. Since we also have
\begin{align*}
  U(t)\( \left\lvert U(-t)\widehat g_a\right\rvert^{p-1}U(-t)\widehat
  g_a\) &=
  \frac{a^2+itp}{a+it}\frac{1}{\(a^2+t^2\)^{n(p-1)/4}}
e^{-\frac{p+it}{a^2+itp}\frac{|x|^2}{2}}\\
&\Tend t 0
  a^{1-n(p-1)/2}e^{-\frac{p}{a^2}\frac{|x|^2}{2}},  
\end{align*}
we check that both terms in
\eqref{eq:poisson} become infinite for $p=1+2/n$, due to a logarithmic
divergence.  

\bibliographystyle{amsplain}
\bibliography{../../carles}

\providecommand{\bysame}{\leavevmode\hbox to3em{\hrulefill}\thinspace}
\providecommand{\MR}{\relax\ifhmode\unskip\space\fi MR }
\providecommand{\MRhref}[2]{%
  \href{http://www.ams.org/mathscinet-getitem?mr=#1}{#2}
}
\providecommand{\href}[2]{#2}
\begin{thebibliography}{1}

\bibitem{COMRL}
R.~Carles and T.~Ozawa, \emph{On the wave operators for the critical nonlinear
  {S}chr\"odinger equation critical nonlinear {S}chr\"odinger equation}, Math.
  Res. Lett. (2007), To appear.

\bibitem{CazCourant}
T.~Cazenave, \emph{Semilinear {S}chr\"odinger equations}, Courant Lecture Notes
  in Mathematics, vol.~10, New York University Courant Institute of
  Mathematical Sciences, New York, 2003.

\bibitem{Ginibre}
J.~Ginibre, \emph{An introduction to nonlinear {S}chr\"odinger equations},
  Nonlinear waves (Sapporo, 1995) (R.~Agemi, Y.~Giga, and T.~Ozawa, eds.),
  GAKUTO International Series, Math. Sciences and Appl., Gakk\={o}tosho, Tokyo,
  1997, pp.~85--133.

\bibitem{KT}
M.~Keel and T.~Tao, \emph{Endpoint {S}trichartz estimates}, Amer. J. Math.
  \textbf{120} (1998), no.~5, 955--980.

\bibitem{Tata1}
D.~Mumford, \emph{Tata lectures on theta. {I}}, Progress in Mathematics,
  vol.~28, Birkh\"auser Boston Inc., Boston, MA, 1983, With the assistance of
  C. Musili, M. Nori, E. Previato and M. Stillman.

\bibitem{NakanishiOzawa}
K.~Nakanishi and T.~Ozawa, \emph{Remarks on scattering for nonlinear
  {S}chr\"odinger equations}, NoDEA Nonlinear Differential Equations Appl.
  \textbf{9} (2002), no.~1, 45--68.

\end{thebibliography}

\end{document}